\documentclass{amsart}
\usepackage{lastpage}
\usepackage{yhmath,verbatim,mathabx}
\usepackage[all]{xy}

\DeclareFontFamily{U}{mathx}{\hyphenchar\font45}
\DeclareFontShape{U}{mathx}{m}{n}{
	<5> <6> <7> <8> <9> <10>
	<10.95> <12> <14.4> <17.28> <20.74> <24.88>
	mathx10
}{}
\DeclareSymbolFont{mathx}{U}{mathx}{m}{n}
\DeclareFontSubstitution{U}{mathx}{m}{n}
\DeclareMathAccent{\widecheck}{0}{mathx}{"71}
\DeclareMathAccent{\wideparen}{0}{mathx}{"75}

\newcommand{\bdis}{\begin{displaymath}}
\newcommand{\edis}{\end{displaymath}}
\newcommand{\be}{\begin{equation}}
\newcommand{\ee}{\end{equation}}
\newcommand{\mbb}{\mathbb}
\newcommand{\mcal}{\mathcal}

\newcommand{\vp}{\varphi}

\newcommand{\vth}{\vartheta}

\newcommand{\zf}{\zeta\left(\frac{1}{2}+it\right)}

\DeclareMathOperator{\dn}{dn}

\DeclareMathOperator{\im}{Im} 

\DeclareMathOperator{\sn}{sn}
\DeclareMathOperator{\cn}{cn}

\theoremstyle{definition}

\theoremstyle{remark}
\newtheorem{remark}[]{Remark}

\newtheorem*{mydef11}{{\bf Theorem 1}}

\newtheorem*{mydef12}{{\bf Theorem 2}}

\newtheorem*{mydef13}{{\bf Theorem 3}}

\newtheorem*{mydef14}{{\bf Theorem 4}}

\newtheorem*{mydef15}{{\bf Theorem 5}} 

\newtheorem*{mydef16}{{\bf Theorem 6}}

\newtheorem*{mydef8}{{\bf Property}}

\numberwithin{equation}{section}

\begin{document}

\title{Jacob's ladders and infinite set of transmutations of asymptotic complete hybrid formula on level curves in Gauss' plane}  

\author{Jan Moser}

\address{Department of Mathematical Analysis and Numerical Mathematics, Comenius University, Mlynska Dolina M105, 842 48 Bratislava, SLOVAKIA}

\email{jan.mozer@fmph.uniba.sk}


\begin{abstract}
In this paper we have obtained new phenomenon lying in the following: every fixed asymptotic complete hybrid formula (we call it as mother formula) generates infinite set of new formulas (transmutations) such that every new formula expresses a close binding between some subset of $\{|\zeta(s)|\}$ and subset of moduli of certain integral and meromorphic functions. \\ 
\begin{center}
	Dedicated to old alchemists
\end{center}
\end{abstract}
\maketitle

\section{Introduction}

\subsection{} 

Let us remind that the following sets of values 
\be \label{1.1} 
\begin{split}
& \left\{ |\zf|^2\right\} \\ 
& \{f_1(t)\}=\{\sin^2t\},\ \{f_2(t)\}=\{\cos^2t\},\ \{f_3(t)\}=\{\cos 2t\}, \\ 
& t\in [\pi L,\pi L+U],\ U\in (0,\pi/4), L\in \mbb{N}
\end{split}
\ee  
generate the following secondary asymptotic complete hybrid formula (see \cite{8}, (3.7), $k_1=k_2=k_3=1$) 
\be \label{1.2} 
\begin{split}
& \left|\zeta(\frac 12+i\alpha_1^{1,1})\right|^2\sin^2(\alpha_0^{1,1})-\\ 
& - \left\{1+\mcal{O}(\frac{\ln\ln L}{\ln L})\right\}\left|\zeta(\frac 12+i\alpha_1^{2,1})\right|^2\cos^2(\alpha_0^{2,1})+\\ 
& + \left|\zeta(\frac 12+i\alpha_1^{3,1})\right|^2\cos(2\alpha_0^{3,1})=0, \\ 
& \forall\- L\geq L_0>0, 
\end{split}
\ee  
($L_0$ is sufficiently big), where 
\be \label{1.3}
\begin{split}
& \alpha_r^{l,1}=\alpha_r(U,\pi L,f_l),\ r=0,1,\ l=1,2,3, \\ 
& \alpha_0^{l,1}\in (\pi L,\pi L+U), \ 
\alpha_1^{l,1}\in (\overset{1}{\wideparen{\pi L}},\overset{1}{\wideparen{\pi L+U}}), 
\end{split}
\ee  
and 
\bdis 
[\overset{1}{\wideparen{\pi L}},\overset{1}{\wideparen{\pi L+U}}]
\edis  
is the first reverse iteration (by means of the Jacob's ladder, see \cite{4}) of the basic segment 
\bdis 
[\pi L,\pi L+U]=[\overset{0}{\wideparen{\pi L}},\overset{0}{\wideparen{\pi L+U}}]. 
\edis 

\begin{remark}
The components of the main $\zeta$-disconnected set (for our case) 
\bdis 
\Delta(\pi L,U,1)=[\pi L,\pi L+U]\cup [\overset{1}{\wideparen{\pi L}},\overset{1}{\wideparen{\pi L+U}}] 
\edis  
are separated each from other by gigantic distance $\rho$, (see \cite{4}, (5.12), comp. \cite{7}, (2.2) -- (2.9))): 
\be \label{1.4} 
\rho\{[\pi L,\pi L+U];[\overset{1}{\wideparen{\pi L}},\overset{1}{\wideparen{\pi L+U}}]\}\sim 
\pi (1-c)\frac{L}{\ln L},\ L\to\infty , 
\ee 
($c$ stands for the Euler's constant).  
\end{remark} 

\subsection{} 

In this paper we obtain, for example, the following transmutations of the mother formula (\ref{1.2}). For every fixed and admissible $U$ and $L$ (see (\ref{1.1}), (\ref{1.2})), the formula (\ref{1.2}) generates: 
\begin{itemize}
	\item[(a)] continuum sets 
	\be \label{1.5} 
	\overset{(1)}{\Omega}_l(U,L), \overset{(2)}{\Omega}_l(U,L),\subset \mbb{C},\ l=1,2,3 
	\ee 
	such that for any elements 
	\bdis 
	s_l^1\in \overset{(1)}{\Omega}_l, s_l^2\in \overset{(2)}{\Omega}_l
	\edis 
	we have the following formula (transmutation of (\ref{1.2})) 
	\be \label{1.6} 
	\begin{split}
	& |\zeta(s_1^1)|^2|\zeta(s_1^2)|^2-\left\{1+\mcal{O}\left(\frac{\ln \ln L}{\ln L}\right)\right\}|\zeta(s_2^1)|^2|\zeta(s_2^2)|^2+|\zeta(s_3^1)|^2|\zeta(s_3^2)|=0, 
	\end{split} 
	\ee 
	\item[(b)] the continuum set 
	\be \label{1.7} 
	\overset{(5)}{\Omega}_l(U,L)\subset\mbb{C}, l=1,2,3 
	\ee 
	such that for any elements 
	\bdis 
	s_l^1\in \overset{(1)}{\Omega}_l, s_l^5\in \overset{(5)}{\Omega}_l
	\edis 
	we have the following formula (transmutation of (\ref{1.2})) 
	\be \label{1.8} 
	\frac{|\zeta(s_1^1)|^2}{|\Gamma(s_1^5)|^2}-\left\{1+\mcal{O}\left(\frac{\ln \ln L}{\ln L}\right)\right\}	\frac{|\zeta(s_2^1)|^2}{|\Gamma(s_2^5)|^2}+	\frac{|\zeta(s_3^1)|^2}{|\Gamma(s_3^5)|}=0, 
	\ee 
	\item[(c)] the continuum set 
	\be \label{1.9} 
	\overset{(7)}{\Omega}_l(U,L,k_l)\subset\mbb{C},\ l=1,2,3,\ [(k_1)^2,(k_2)^2,(k_3)^2]\in (0,1)^3
	\ee 
	(for every fixed and admissible $k_1,k_2,k_3$) such that for any fixed elements 
	\bdis 
	s_l^1\in \overset{(1)}{\Omega}_l, s_l^7\in \overset{(5)}{\Omega}_l
	\edis 
	we have the following transmutation of (\ref{1.2}) 
	\be \label{1.10} 
	\begin{split}
	& |\zeta(s_1^1)|^2|\sn(s_1^7,k_1)|^2-\\ 
	& -\left\{1+\mcal{O}\left(\frac{\ln \ln L}{\ln L}\right)\right\}|\zeta(s_2^1)|^2|\cn(s_2^7,k_2)|^2+|\zeta(s_3^1)|^2|\dn(s_3^7,k_3)|=0, 
	\end{split} 
	\ee 
	(for Jacobi's elliptic functions of moduli $k_1,k_2,k_3$). 
	\end{itemize} 

\subsection{} 

Now we give the following remarks. 

\begin{remark}
The formula (\ref{1.1}) gives the example of the continuum set of transmutations of the mother formula (\ref{1.2}). 
\end{remark} 

\begin{remark}
Let us remind that the mother formula is written as orthogonality condition for corresponding vectors. We see that the transmutations (\ref{1.6}), (\ref{1.8}) as well as (\ref{1.10}) inherit this property. 
\end{remark} 

\begin{remark}
The method presented in this paper is general one in the following sense: 
\begin{itemize}
	\item[(a)] single secondary asymptotic complete hybrid formula (\ref{1.2}) generates the infinite set of its transmutations, 
	\item[(b)] there is an infinite set of asymptotic complete hybrid formulas, (see \cite{5} -- \cite{8} concerning evolution from $\zeta$-factorization formula to complete hybrid formula). 
\end{itemize}
\end{remark}

\begin{remark}
Let us remind that the general feature of all transmutations is a new phenomenon that lies in close binding between some subsets of $|\zeta(s)|$ and subsets of moduli of other integral and meromorphic functions. 
\end{remark} 

\begin{remark}
This paper is based also on new notions and methods in the theory of the Riemann zeta-function we have introduced in our series of 48 papers concerning Jacob's ladders. These can be found in arXiv [math.CA] starting with the paper \cite{2}. 
\end{remark} 

\section{Jacob's ladder as basis of secondary complete hybrid formula} 

\subsection{} 

Let us remind that the Jacob's ladder 
\bdis 
\vp_1(t)=\frac 12\vp(t)
\edis  
was introduced in \cite{2}, (see also \cite{3}), where the function $\vp_1(t)$ is arbitrary solution of the non-linear integral equation (also introduced in \cite{2}) 
\bdis 
\int_0^{\mu[x(T)]}Z^2(t)e^{-\frac{2}{x(T)}t}{\rm d}t=\int_0^T Z^2(t){\rm d}t, 
\edis 
of course, 
\bdis 
\begin{split}
& Z(t)=e^{i\vth(t)}\zf,\\ 
& \vth(t)=-\frac t2\ln\pi+\im\left\{\ln\Gamma\left(\frac 14+i\frac t2\right)\right\}, 
\end{split}
\edis 
where each admissible function $\mu(y)$ generates a solution 
\bdis 
y=\vp(T;\mu)=\vp(T);\ \mu(y)>7y\ln y. 
\edis 
We call the function $\vp_1(T)$ the Jacob's ladder as an analogue of the Jacob's dream in Chumash, Bereishis, 28:12. 

\begin{remark}
By making use of the Jacob's ladders we have shown (see \cite{2}) that the classical Hardy-Littlewood integral (1918) 
\bdis 
\int_0^T\left|\zf\right|{\rm d}t
\edis 
has - in addition to previously known Hardy-Littlewood expression (and other similar) possessing an unbounded error term at $T\to\infty$ - the following infinite set of almost exact representations 
\bdis 
\int_0^T\left|\zf\right|{\rm d}t=\vp_1(T)\ln\{\vp_1(T)\}+(c-\ln 2\pi)\vp_1(T)+c_0+\mcal{O}\left(\frac{\ln T}{T}\right) 
\edis  
as 
\bdis 
T\to\infty , 
\edis  
where $c$ is the Euler's constant and $c_0$ is the constant from the Titschmarsh-Kober-Atkinson formula. 
\end{remark} 

\subsection{} 

Next, we have obtained (see \cite{2}, (6.2)) the following formula 
\be \label{2.1} 
T-\vp_1(T)\sim (1-c)\pi (T);\ \pi(T)\sim \frac{T}{\ln T},\ T\to\infty , 
\ee  
where $\pi(T)$ is the prime-counting function. 

\begin{remark}
Consequently, the Jacob's ladder $\vp_1(T)$ can be viewed by our formula (\ref{2.1}) as an asymptotic complementary function to the function 
\bdis 
(1-c)\pi(T) 
\edis  
in the following sense 
\bdis 
\vp_1(T)+(1-c)\pi(T)\sim T,\ T\to\infty. 
\edis 
\end{remark} 

Since Jacob's ladder is exactly increasing function we have the reversely iterated sequence 
\bdis 
\{\overset{k}{T}\}_{k=1}^{k_0},\ k_0\in\mbb{N}:\ \vp_1(\overset{k}{T})=\overset{k-1}{T},\ \overset{0}{T}=T>T_0>0,\ \overset{k}{T}=\vp_1^{-k}(T), 
\edis  
where $T_0$ is a sufficiently big (we fix $k_0$). Next we have (see \cite{4}, (2.5) -- (2.7), (5.1) -- (5.13)) the following basic property generated by the Jacob's ladder. 

\begin{mydef8}
For every segment 
\bdis  
[T,T+U],\ U=o\left(\frac{T}{\ln T}\right),\ T\to\infty 
\edis 
there is the following class of disconnected sets 
\be \label{2.2} 
\Delta(T,U,k)=\bigcup_{r=0}^k [\overset{r}{\wideparen{T}},\overset{r}{\wideparen{T+U}}],\ 1\leq k\leq k_0
\ee 
generated by the Jacob's ladder $\vp_1(t)$. 
\end{mydef8} 

\begin{remark}
Disconnected set (\ref{2.2}) has the following properties (see \cite{4}, (2.5) -- (2.7)): 
\begin{itemize}
	\item[(a)] lengths of its components are given by 
	\be \label{2.3} 
	|[\overset{r}{\wideparen{T}},\overset{r}{\wideparen{T+U}}]|\sim o\left(\frac{T}{\ln T}\right),\ T\to\infty, 
	\ee 
	\item[(b)] lengths of adjacent segments are given by 
	\be \label{2.4} 
	|[\overset{r-1}{\wideparen{T+U}},\overset{r}{\wideparen{T}}]|\sim (1-c)\frac{T}{\ln T},\ T\to\infty, 
	\ee 
	i.e. for distance $\rho$ of two consecutive components we have 
	\be \label{2.5} 
	\rho\{[\overset{r-1}{\wideparen{T}},\overset{r-1}{\wideparen{T+U}}];[\overset{r}{\wideparen{T}},\overset{r}{\wideparen{T+U}}]\}\sim (1-c)\frac{T}{\ln T},\ T\to\infty, \ r=1,\dots,k. 
	\ee 
\end{itemize}
\end{remark} 

\begin{remark}
Formula (\ref{1.4}) follows from (\ref{2.5}) for $r=1$ and $T=\pi L$. 
\end{remark} 

\begin{remark}
Asymptotic behavior of the components of the disconnected set (\ref{2.2}) is as follows: at $T\to\infty$ these components recede unboundedly each from other and all together are receding to infinity. Hence, at $T\to\infty$ the set of the components of (\ref{2.2}) behaves as one-dimensional Friedman-Hubble universe. 
\end{remark} 

\subsection{} 

Finally, let us remind that the functions (see (\ref{1.2}), (\ref{1.3})) 
\bdis 
\alpha_r^{l,1}=\alpha_r(U,\pi L,f_l),\ r=0,1,\ l=1,2,3 
\edis  
are generated by the Jacob's ladder $\vp_1(T)$ and its iterations as follows (see \cite{4}, (6.3), \cite{5}, (1.14), \cite{7}, (3.6)) 
\be \label{2.6}  
\begin{split}
& \alpha_r^{l,1}=\vp_1^{1-r}(d_l),\ d_l=d(U,\pi L,f_l),\ r=0,1,\ l=1,2,3, \\ 
& \alpha_r^{l,1}\in (\overset{r}{\wideparen{\pi L}},\overset{r}{\wideparen{\pi L+U}}),\ r=0,1. 
\end{split}
\ee  

\section{Transmutation of type $[\zeta,\zeta]$} 

\subsection{} 

Let us remind that the loci 
\be \label{3.1} 
|\zeta(s)|=c,\ c>0,\ s\in\mbb{C}\setminus \{1\} 
\ee  
is the level curve, (comp. \cite{9}, pp. 120 -- 122). Of course, the topological structure of the level curve is (in general) very complicated. For example, the level curve 
\bdis 
|\zeta(s)|=0.8 
\edis  
contains three ovals (apart from another) situated in the neighborhood of the first three roots of the equation 
\bdis 
\zf=0, 
\edis  
i.e. this one represents the disconnected set, (see \cite{1}, p. 371, FIG. 149). 

\subsection{} 

First, we define for every admissible and fixed $U,L$ (see (\ref{1.1}), (\ref{1.2})) the following level curves: 
\be \label{3.2} 
\overset{(1)}{\Omega}_l(U,L),\ l=1,2,3 
\ee  
as the loci 
\be \label{3.3} 
|\zeta(s_l^1)|=\left|\zeta\left( \frac 12+i\alpha_1^{l,1}\right)\right|, 
\ee 
where 
\be \label{3.4} 
s_l^1=s_l^1(U,L)\in \overset{(1)}{\Omega}_l\subset\mbb{C},\ l=1,2,3, 
\ee  
(of course, the sets $\overset{(1)}{\Omega}_l$ are continuum ones). 

\subsection{} 

Secondly, we use the loci (\ref{3.1}) for the following values 
\be \label{3.5} 
c_1=|\sin \alpha_0^{1,1}|,\ c_2=|\cos \alpha_0^{2,1}|,\ c_3=\cos 2\alpha_0^{3,1}, 
\ee  
see (\ref{1.2}) where, of course 
\bdis 
\cos (2\alpha_0^{3,1})=\cos[2(\alpha_0^{3,1}-\pi L)]>0. 
\edis 
Consequently, we define for every admissible and fixed $U,L$ the following three level curves 
\be \label{3.6} 
\overset{(2)}{\Omega}_l(U,L),\ l=1,2,3 
\ee  
as the loci 
\be \label{3.7} 
\begin{split}
& |\zeta(s_1^2)|=|\sin\alpha_0^{1,1}|,\ s_1^2=s_1^2(U,L)\in \overset{(2)}{\Omega}_1, \\ 
& |\zeta(s_2^2)|=|\cos\alpha_0^{2,1}|,\ s_2^2=s_2^2(U,L)\in \overset{(2)}{\Omega}_2, \\ 
& |\zeta(s_3^2)|=\cos 2\alpha_0^{3,1},\ s_3^2=s_3^2(U,L)\in \overset{(2)}{\Omega}_3. 
\end{split}
\ee 

\subsection{} 

Finally, construction of the sets $\overset{(1)}{\Omega}_l,\overset{(2)}{\Omega}_l$ together with application of the substitutions (\ref{2.3}), (\ref{3.7}) into (\ref{1.2}) completes the proof of the following. 

\begin{mydef11}
For every fixed and admissible $U,L$ (see (\ref{1.1}), (\ref{1.2})) there are continuum sets 
\be \label{3.8} 
\overset{(1)}{\Omega}_l(U,L),\ \overset{(2)}{\Omega}_l(U,L)\subset\mbb{C},\ l=1,2,3 
\ee  
such that for every elements 
\be \label{3.9} 
s_l^1\in \overset{(1)}{\Omega}_l,\ s_l^2\in \overset{(2)}{\Omega}_l 
\ee  
we have the following formula (transmutation of (\ref{1.2})) 
\be \label{3.10} 
\begin{split}
& |\zeta(s_1^1)|^2|\zeta(s_1^2)|^2-\left\{1+\mcal{O}\left(\frac{\ln \ln L}{\ln L}\right)\right\}|\zeta(s_2^1)|^2|\zeta(s_2^2)|^2+|\zeta(s_3^1)|^2|\zeta(s_3^2)|=0. 
\end{split}
\ee 
\end{mydef11} 

\section{Transmutations of types $[\zeta,\cos]$, $[\zeta,(s)^{n_1},(s)^{n_2},(s)^{n_3}]$} 

\subsection{} 

We use the loci 
\be \label{4.1} 
|\cos s|=c,\ c>0,\ s\in\mbb{C} 
\ee  
for the three values (\ref{3.5}). Consequently, we define for every admissible and fixed $U,L$ the following level curves 
\be \label{4.2} 
\overset{(3)}{\Omega}_l,\ l=1,2,3 
\ee  
as the loci 
\be \label{4.3} 
\begin{split}
	& |\cos s_1^3|=|\sin\alpha_0^{1,1}|,\ s_1^3=s_1^3(U,L)\in \overset{(3)}{\Omega}_1, \\ 
	& |\cos s_2^3|=|\cos\alpha_0^{2,1}|,\ s_2^3=s_2^3(U,L)\in \overset{(3)}{\Omega}_2, \\ 
	& |\cos s_3^3|=\cos 2\alpha_0^{3,1},\ s_3^3=s_3^3(U,L)\in \overset{(3)}{\Omega}_3. 
\end{split}
\ee 
Now we have the following. 

\begin{mydef12}
For every fixed and admissible $U,L$ there are the sets 
\be \label{4.4} 
\overset{(1)}{\Omega}_l(U,L),\ \overset{(3)}{\Omega}_l(U,L)\subset\mbb{C},\ l=1,2,3 
\ee  
such that we have the following formula (transmutation of (\ref{1.2})) 
\be \label{4.5} 
\begin{split}
	& |\zeta(s_1^1)|^2|\cos s_1^3|^2-\left\{1+\mcal{O}\left(\frac{\ln \ln L}{\ln L}\right)\right\}|\zeta(s_2^1)|^2|\cos s_2^3|^2+|\zeta(s_3^1)|^2|\cos s_3^3|=0. 
\end{split}
\ee 
\end{mydef12} 

\subsection{} 

For the functions 
\be \label{4.6} 
(s)^{n_1},\ (s)^{n_2},\ (s)^{n_3},\ (n_1,n_2,n_3)\in\mbb{N}^3 
\ee  
we define the level curves 
\be \label{4.7} 
\overset{(4)}{\Omega}_l(U,L,n_l),\ l=1,2,3 
\ee  
as the loci 
\be \label{4.8} 
\begin{split} 
& |s_1^4|^{n_1}=|\sin \alpha_0^{1,1}|, \\ 
& |s_2^4|^{n_2}=|\cos \alpha_0^{2,1}|, \\ 
& |s_3^4|^{n_3}=\cos 2\alpha_0^{3,1} 
\end{split} 
\ee  
for every fixed and admissible $U,L,n_1,n_2,n_3$. Now we have (see (\ref{1.2}), (\ref{3.4}), (\ref{4.8})) the following. 

\begin{mydef13}
For every fixed and admissible $U,L,n_1,n_2,n_3$ there are the sets 
\be \label{4.9} 
\overset{(1)}{\Omega}_l(U,L), \overset{(4)}{\Omega}_l(U,L,n_l)\subset\mbb{C},\ l=1,2,3 
\ee  
such that we have the following formula (transmutation of (\ref{1.2})) 
\be \label{4.10} 
\begin{split}
	& |\zeta(s_1^1)|^2|s_1^4|^{2n_1}-\left\{1+\mcal{O}\left(\frac{\ln \ln L}{\ln L}\right)\right\}|\zeta(s_2^1)|^2|s_2^4|^{2n_2}+|\zeta(s_3^1)|^2|s_3^4|^{n_3}=0. 
\end{split}
\ee 
\end{mydef13} 

\begin{remark}
The formula (\ref{4.10}) gives an example of infinite set of transmutations of the mother formula (\ref{1.2}). 
\end{remark} 

\section{Transmutations of types $[\zeta,\frac{1}{\Gamma}]$, $[\zeta,J_{p_1},J_{p_2},J_{p_3}]$, $[\zeta,\sn,\cn,\dn]$} 

\subsection{} 

For the integral function 
\be \label{5.1} 
\frac{1}{\Gamma(s)},\ s\in\mbb{C} 
\ee  
we define the level curves 
\be \label{5.2} 
\overset{(5)}{\Omega}_l(U,L),\ l=1,2,3 
\ee  
as the loci 
\be \label{5.3} 
\begin{split}
& \frac{1}{|\Gamma(s_1^5)|}=|\sin \alpha_0^{1,1}|, \\ 
& \frac{1}{|\Gamma(s_2^5)|}=|\cos \alpha_0^{2,1}|, \\ 
& \frac{1}{|\Gamma(s_3^5)|}=\cos 2\alpha_0^{3,1}, 
\end{split}
\ee 
for every fixed and admissible $U,L$. Now we have (see (\ref{1.2}), (\ref{3.4}), (\ref{5.3})) the following 

\begin{mydef14}
For every fixed and admissible $U,L$ there are sets 
\be \label{5.4} 
\overset{(1)}{\Omega}_l(U,L),\ \overset{(5)}{\Omega}_l(U,L)\subset\mbb{C},\ l=1,2,3 
\ee  
such that we have the following formula (transmutation of (\ref{1.2})) 
\be \label{5.5} 
\frac{|\zeta(s_1^1)|^2}{|\Gamma(s_1^5)|^2}-\left\{1+\mcal{O}\left(\frac{\ln \ln L}{\ln L}\right)\right\}\frac{|\zeta(s_2^1)|^2}{|\Gamma(s_2^5)|^2}+\frac{|\zeta(s_3^1)|^2}{|\Gamma(s_3^5)|}=0. 
\ee 
\end{mydef14} 

\subsection{} 

For the Bessel's functions 
\be \label{5.6} 
J_{p_1}(s), J_{p_2}(s), J_{p_3}(s), s\in\mbb{C}, (p_1,p_2,p_3)\in\mbb{Z}^3
\ee 
we define the level curves 
\be \label{5.7} 
\overset{(6)}{\Omega}_l(U,L,p_l),\ l=1,2,3 
\ee  
as the loci 
\be \label{5.8} 
\begin{split}
& |J_{p_1}(s_1^6)|=|\sin \alpha_0^{1,1}|, \\ 
& |J_{p_2}(s_2^6)|=|\cos \alpha_0^{2,1}|, \\ 
& |J_{p_3}(s_3^6)|=\cos 2\alpha_0^{3,1}, 
\end{split}
\ee  
for every fixed and admissible $U,L,p_1,p_2,p_3$. Now we have (see (\ref{1.2}), (\ref{3.4}), (\ref{5.8})) the following. 

\begin{mydef15}
For every fixed and admissible $U,L,p_1,p_2,p_3$ there are sets 
\be \label{5.9} 
\overset{(1)}{\Omega}_l(U,L), \overset{(6)}{\Omega}_l(U,L,p_l)\subset\mbb{C}, l=1,2,3 
\ee 
such that we have the following formula (transmutation of (\ref{1.2})) 
\be \label{5.10} 
\begin{split}
	& |\zeta(s_1^1)|^2|J_{p_1}(s_1^6)|-\left\{1+\mcal{O}\left(\frac{\ln \ln L}{\ln L}\right)\right\}|\zeta(s_2^1)|^2|J_{p_2}(s_2^6)|^{2}+ \\ 
	& + |\zeta(s_3^1)|^2|J_{p_3}(s_3^6)|=0. 
\end{split}
\ee 
\end{mydef15} 

\begin{remark}
Here we have the infinite (enumerable) set of transmutations of transmutations of the mother formula (\ref{1.2}), comp. Remark 12. 
\end{remark}

\subsection{} 

For the Jacobi's elliptic functions (their relief-surfaces with $k=0.8$ see in \cite{1}, pp. 202, 203) 
\be \label{5.11} 
\sn(s,k_1),\ \cn(s,k_2),\ \dn(s,k_3),\quad [(k_1)^2,(k_2)^2,(k_3)^2]\in (0,1)^3
\ee  
($k_1,k_2,k_3$ stand for the moduli) we define the level curves 
\be \label{5.12} 
\overset{(7)}{\Omega}_l(U,L,k_l), \ l=1,2,3 
\ee  
as the loci 
\be \label{5.13} 
\begin{split}
& |\sn(s_1^7,k_1)|=|\sin \alpha_0^{1,1}|, \\ 
& |\cn(s_2^7,k_2)|=|\cos \alpha_0^{2,1}|, \\ 
& |\dn(s_3^7,k_3)|=\cos 2\alpha_0^{3,1}
\end{split}
\ee 
for every fixed and admissible $U,L,k_1,k_2,k_3$. Now we have (see (\ref{1.2}), (\ref{3.4}), (\ref{5.13})) the following. 

\begin{mydef16} 
For every fixed and admissible $U,L,k_1,k_2,k_3$ there are sets 
\be \label{5.14} 
\overset{(1)}{\Omega}_l(U,L), \overset{(7)}{\Omega}_l(U,L,k_l)\subset\mbb{C}, l=1,2,3 
\ee 
such that we have the following formula (transmutation of (\ref{1.2})) 
\be \label{5.15} 
\begin{split}
	& |\zeta(s_1^1)|^2|\sn(s_1^7,k_1)|-\left\{1+\mcal{O}\left(\frac{\ln \ln L}{\ln L}\right)\right\}|\zeta(s_2^1)|^2|\cn(s_2^7,k_2)|^{2}+ \\ 
	& + |\zeta(s_3^1)|^2|\dn(s_3^7,k_3)|=0. 
\end{split}
\ee 
\end{mydef16}  

\begin{remark}
About the potency of the set of transmutations in (\ref{5.15}) see Remark 2. 
\end{remark}

I would like to thank Michal Demetrian for his moral support of my study of Jacob's ladders.

\end{document}